\documentclass[11pt]{article}
\usepackage{amsmath,amssymb,amsfonts}
\usepackage{enumerate}
\usepackage[latin1]{inputenc}
\allowdisplaybreaks
\usepackage{color}
\usepackage[colorlinks]{hyperref}

\usepackage{cancel}

\author{Fazia {\sc Bedouhene}\thanks{Mouloud Mammeri University of Tizi-Ouzou,
Laboratoire de Math\'ematiques Pures et Appliqu\'ees, Tizi-Ouzou, Algeria
E-Mail: fbedouhene@yahoo.fr},
Youcef {\sc Ibaouene}\footnotemark[1] \thanks{Normandie Univ, Laboratoire Rapha\"el Salem,
UMR CNRS 6085, Rouen, France
E-Mail: youcef.ibaouene@etu.univ-rouen.fr},
Omar {\sc Mellah}\thanks{Mouloud Mammeri University of Tizi-Ouzou,
Laboratoire de Math\'ematiques Pures et Appliqu\'ees, Tizi-Ouzou, Algeria
E-Mail: omellah@yahoo.fr},
\and
Paul {\sc Raynaud de Fitte} \thanks{Normandie Univ, Laboratoire Rapha\"el Salem,
UMR CNRS 6085, Rouen, France E-Mail: prf@univ-rouen.fr}}
\title{Weyl almost periodic solutions to abstract linear and semilinear equations with Weyl almost periodic
coefficients }

\newtheorem{theo}{Theorem}[section]
\newtheorem{prop}[theo]{Proposition}

\newtheorem{lemma}[theo]{Lemma}
\newtheorem{definition}[theo]{Definition}

\newtheorem{remark}[theo]{Remark}

\newtheorem{example}[theo]{Example}

\renewcommand\epsilon{\varepsilon}

\newcommand\N{\mathbb{N}}
\newcommand\R{\mathbb{R}}

\DeclareMathOperator*{\argmax}{arg\,max}

\newcommand\un[1]{\,\rlap{{1}}\kern.22em \mbox{l}_{#1}} 

\newcommand\proof{\noindent {\bf Proof}\ }

\renewcommand\square{\fbox{\rule{0em}{.3em}\rule{.3em}{0em}} \qquad}
\newcommand\finpr{\hfill$\square\qquad$\medskip\par}

\newcommand\CCO[1]{\left( #1 \right)}
\newcommand\norm[1]{\left\Vert #1 \right\Vert}
\newcommand\abs[1]{\left\vert #1 \right\vert}
\newcommand\accol[1]{\left\{#1\right\}}

\newcommand\loc{\mathrm{loc}}

\newcommand\CB{\mbox{\rm CB}}

\newcommand\StAP[2][]{\St^{#2}_{#1}\mbox{\rm{AP}}} 

\newcommand\St{\mathbb{S}} 
\newcommand\We{\mathbb{W}} 

\newcommand\WeAP[2][]{\We^{#2}_{#1}\mbox{\rm{AP}}} 

\newcommand\espX{{\mathbb X}}

\newcommand{\verti}[1]{{\left\vert\kern-0.25ex\left\vert\kern-0.25ex\left\vert #1
    \right\vert\kern-0.25ex\right\vert\kern-0.25ex\right\vert}}

\begin{document}

\maketitle

\begin{abstract}
In this work,   we study the existence  and uniqueness of  bounded Weyl almost periodic solution to the abstract differential equation
$u'(t) = Au(t) + f(t)$, $t\in\R$, in a Banach space $\espX$, where $A:D\left( A\right) \subset \espX\rightarrow \espX$ is a linear operator  (unbounded) which generates an exponentially stable $C_{0}$-semigroup on $\espX$ and $f :\R \rightarrow \espX $ is a Weyl almost periodic function. We also investigate the nonautonomous case.
\end{abstract}



\section{Introduction}

In 1927, Hermann Weyl \cite{Weyl_1927} introduced a generalization of Bohr and Stepanov almost periodic functions. Since then, generalized almost periodic functions bear his name and numerous contributions to theory of almost periodic functions have been brought, see for instance \cite{Andres_Bersani_Grande_06_Hierarchy, Andres_2001, besicovitch1931almost, besicovitoh1932almost, danilov2006weyl, Levitan_1966,Ursell1931_Parseval}.  Weyl-almost periodic functions are very important and have several applications in dynamical systems, in particular in symbolic dynamics \cite{iwanik1988weyl}, where an important class of Weyl almost periodic trajectories are the regular Toeplitz sequences of Jacobs and Keane \cite{jacobs19690}.

In recent years, many authors have been interested in the study of existence and uniqueness of almost periodic solutions to different types of differential equations with Stepanov almost periodic coefficients in both deterministic and stochastic cases. We can mention, in the deterministic case, the work by Andres and Pennequin \cite{andres2012stepanov,andres2012nonexistence}, Long and Ding \cite{long2011composition}, Ding et al. \cite{Ding_Long_Nguerekata2011}, Hu and Mingarelli  \cite{hu_2008_bochner}, Henriquez \cite{henriquez1990stepanov}, Rao \cite{rao1975stepanov}, Zaidman \cite{zaidman1971existence}. Particularly, in \cite{Ding_Long_Nguerekata2011}, \cite{long2011composition}, and  \cite{Ait_Dads_Ezzinbi_16_MMAS}, the authors show the existence and uniqueness of almost periodic solution for an  abstract semilinear evolution equation with Stepanov almost periodic coefficients. They prove the existence and uniqueness of a Bohr almost periodic mild solution.
Nevertheless, in the work \cite{andres2012nonexistence}, Andres and Pennequin  proved the nonexistence of purely Stepanov almost periodic solutions of ordinary differential equations in uniformly convex Banach spaces.
Until now, to our knowledge, there are only few works dedicated to the study of existence and uniqueness of Weyl almost periodic solutions to some differential equations. The first investigation in this direction is due to Lenka Radova's \cite{Radova_04}, but her study was limited to the scalar equation  $x' = ax + f(t)$, where $a \in \R$ and $f:\R\rightarrow\R$ is an essentially bounded Weyl almost periodic function. The second one is due to Andres et al. \cite{andres_Bersani_Radova2006almost}, and can be seen as an extension of Radova's work. Precisely, the authors investigate almost periodic solutions in various senses (Stepanov, Weyl, Besicovitch) of higher-order scalar differential equations with a nonlinear Lipschitz restoring term and a bounded Stepanov, Weyl or Besicovitch almost periodic forcing term. One can also mention the recent paper by Kostic \cite{kostic_generalized}, where Radova's study was addressed in the context of linear evolution equations with bounded and (asymptotically) Weyl-almost periodic coefficients.

Motivated by the previous papers, the aim of this paper is to generalize the previous investigations to the case of semilinear evolution equations, without any boundedness restriction on the forcing term. More precisely, we consider the following abstract semilinear evolution equation:
\begin{equation}\label{eq:EDO semi-lineaire}
u^{\prime }\left( t\right) = Au(t) +f( t,u(t))\quad\quad\forall t\in\R,
\end{equation}
where $A:D\left( A\right) \subset \espX\rightarrow \espX$ is a linear operator (possibly unbounded) which generates an exponentially stable $C_{0}$-semigroup on a Banach space $\espX$, and $(t,x)\mapsto f(t,x)$  is a parametric Weyl almost periodic function of degree $p\geq 1$. We show that, under some conditions,  \eqref{eq:EDO semi-lineaire} has a unique solution which is bounded and Weyl almost periodic. Note that, contrarily to the Stepanov case, the solution can be purely Weyl almost periodic.

 The major problem encountered in our study is that the space of Weyl almost periodic functions endowed with the Weyl semi-norm is not complete \cite{Andres_2001}, which makes the approach and classical tools of functional analysis inapplicable, especially the fixed-point theorem in Banach spaces. To overcome this difficulty, we have opted for the method used by Kamenskii et al \cite{Kamenskii_Mikhail_and_Mellah_Omar_and_Raynaud_de_Fitte_15}, which consists in showing "manually" that the unique bounded mild solution to \eqref{eq:EDO semi-lineaire}, obtained using the fixed-point theorem,  is Weyl almost periodic. 

This work is organized as follows: In the second section, we recall some definitions and results related to Stepanov almost periodic and Weyl almost periodic functions. In section 3, we present three essential parts. In the first part, based on a compactness property established by Danilov \cite{danilov2006weyl}, we give a new superposition result in the space of Weyl almost periodic functions. In the second part, after providing examples of bounded and purely Weyl almost periodic solutions, we tackle the problem of existence and uniqueness of a bounded Weyl almost periodic solution to an abstract linear differential equation. Finally, in the third part, inspired from \cite{Kamenskii_Mikhail_and_Mellah_Omar_and_Raynaud_de_Fitte_15}, we prove that equation \eqref{eq:EDO semi-lineaire} has a unique bounded mild  solution which is Weyl almost periodic.
\section{Notations and Preliminaries}
 First, we give some basic definitions and results  on Stepanov almost periodic functions and Weyl almost periodic functions (for more details, see \cite{Andres_Bersani_Grande_06_Hierarchy, besicovitch1931almost, besicovitoh1932almost}).

We denote by $\R$ the set of real numbers and by $\espX$ a Banach space endowed with the norm $\norm{.}$. 

We denote also by $\CB{(\R, \espX)}$ the Banach space of continuous and bounded functions from $\R$ to $\espX$, endowed with the norm
\begin{equation*}
\norm{u}_{\infty} = \sup_{t\in\R}\norm{u(t)}.
\end{equation*}
\subsection{Stepanov and Weyl norm}
Given $l> 0$, the Stepanov norm associated with $l$ of a locally integrable function $f: \R\rightarrow \espX$, i.e $f\in L^{p}_{\loc}(\R, \espX)$ $(p\geq 1)$, is defined by
\begin{equation*}
\norm{f}_{\St^{p}_{l}}= \sup_{\xi\in \R}\CCO{\frac{1}{l}\int_{\xi}^{\xi+l}\norm{f(t)}^{p}dt}^{\frac{1}{p}}.
\end{equation*}
The Weyl norm is defined by
\begin{equation*}
\norm{f}_{\We^{p}}= \lim_{l\rightarrow+\infty}\norm{f}_{\St^{p}_{l}}.
\end{equation*}
The limit always exists for $l\rightarrow +\infty$ \cite[p.72-73]{besicovitoh1932almost}.
The distance induced by $\norm{.}_{\St^{p}_{l}}$ (respectively $\norm{f}_{\We^{p}}$)  is denoted by $D_{\St_{l}^{p}}$ (respectively $D_{\We^{p}}$).
\subsection{Almost periodicity in Stepanov and Weyl senses}
 Let us recall some definitions of Stepanov and Weyl almost periodic functions. Recall that a set $T$ of real numbers is {\em relatively dense} if there exists a real number $k>0$, such that $T\cap[a, a+k]\neq\emptyset$, for all $a$ in $\R$.
 \begin{definition}{\em(\cite[p.76-77]{Amerio_1971}, \cite{Andres_Bersani_Grande_06_Hierarchy,  Andres_2001}, \cite[p.77]{besicovitoh1932almost},  \cite[p.173]{Corduneanu}, \cite[p.188]{Levitan_1966}) A  function  $f\in L_{\loc}^{p}(\R, \espX)$ is said to be {\em almost periodic in the sense of Stepanov} (or {\em $p$-Stepanov almost periodic} if we want to highlight the degree $p$) if, for every  $\epsilon>0$, the set
\begin{equation*}
T_{\St^{p}}(\epsilon,f)=\left\{\tau\in\R;\sup_{\xi\in\R}\CCO{\int_{\xi}^{\xi+1}\norm{f(t+\tau)- f(t)}^{p}dt}^{\frac{1}{p}}<\epsilon\right\}
\end{equation*}
is relatively dense. The number $\tau$ is called an {\em $\epsilon$-Stepanov almost-period} (or {\em Stepanov $\epsilon$-translation number} of $f$).
We denote the set of all such functions by $\StAP{p}(\R, \espX)$.}
\end{definition}
\begin{definition}  {\em (\cite{Andres_Bersani_Grande_06_Hierarchy,  Andres_2001},  \cite[p.78]{besicovitoh1932almost}, \cite[p.190]{Levitan_1966}) A  function  $f\in L_{\loc}^{p}(\R, \espX)$ is called {\em Weyl almost periodic} if, for every  $\epsilon>0$,  there exist $l = l(\epsilon) >0$ such that, the set
\begin{equation}\label{for: definition equi-Weyl}
T_{\St^{p}_{l}}(\epsilon, f)=\left\{\tau\in\R;\sup_{\xi\in\R}\CCO{\frac{1}{l}\int_{\xi}^{\xi+l}\norm{f(t+\tau)- f(t)}^{p}dt}^{\frac{1}{p}}<\epsilon\right\}
\end{equation}
is relatively dense. This means that for every $\epsilon>0$, there exists $l=l(\epsilon)>0$ such that $f\in \StAP[l]{p}(\R, \espX)$. The number $\tau$ is called a {\em Weyl $\epsilon$-translation number} of $f$.  We denote the set of all such functions by $\WeAP{p}(\R,\espX)$.}
\end{definition}

 There is, however, another interesting definition given by Ursell \cite{Ursell1931_Parseval}:
 \begin{definition}[\cite{Ursell1931_Parseval}]\label{def:Ursell}
 {\em A  function  $f\in L_{\loc}^{p}(\R, \espX)$ is said
to be {\em Weyl *-almost periodic} (or W.a.p  following Ursell's notations) if, for every  $\epsilon>0$,  there exists $l = l(\epsilon)$ such that the set
\begin{equation*}
T^{*}_{\St^{p}_{l}}(\epsilon,f)= \left\{\tau\in\R;\CCO{\frac{1}{l}\int_{0}^{l}\norm{f(t + \tau)-  f(t)}^{p}dt}^{\frac{1}{p}}< \epsilon\right\}
\end{equation*}
is relatively dense.  We denote the set of all such functions by $\WeAP[*]{p}(\R,\espX)$.}
\end{definition}
Ursell \cite{Ursell1931_Parseval} has shown that
\begin{theo}\label{Theorem: Ursell}
{\em  The classes $\WeAP{p}(\R,\espX)$ and $\WeAP[*]{p}(\R,\espX)$ coincide.}
 \end{theo}
  This identification is very important when establishing our main result.
 \remark {\em
\begin{enumerate}
  \item Some authors name the Weyl almost periodicity by equi-Weyl almost periodicity (e.g~ \cite{Andres_Bersani_Grande_06_Hierarchy,  Andres_2001,
      Radova_04}), distinguishing between the above definition and the following one introduced by Kovanko \cite{kovanko1944compacite}:
A  function  $f\in L_{\loc}^{p}(\R, \espX)$ is said to be {\em almost periodic in the sense of Weyl} according to \cite{Andres_Bersani_Grande_06_Hierarchy,  Andres_2001, Radova_04} if, for every  $\epsilon>0$, the set
\begin{equation*}
T_{\We^{p}}(\epsilon, f)=\left\{\tau\in\R;\lim_{l\rightarrow+\infty}\sup_{\xi\in\R}\CCO{\frac{1}{l}\int_{\xi}^{\xi+l}\norm{f(t+\tau)- f(t)}^{p}dt}^{\frac{1}{p}}<\epsilon\right\}
\end{equation*}
is relatively dense. We denote the set of all such functions by $\We^{p}(\R, \espX)$.

Of course the space $\WeAP{p}(\R,\espX)$ is an intermediate space between $\StAP{p}(\R, \espX)$ and $\We^{p}(\R,\espX)$, i.e. we have
\begin{equation*}
 \StAP{p}(\R, \espX) \subset \WeAP{p}(\R,\espX)\subset \We^{p}(\R,\espX).
\end{equation*}
 According to \cite{Andres_Bersani_Grande_06_Hierarchy, Andres_2001}, the previous inclusions are strict.
  \item Every Weyl almost periodic function of degree $p$ is $\St^{p}$-bounded and also $\We^{p}$-bounded and also  (see \cite{Andres_Bersani_Grande_06_Hierarchy}).
 Actually, the spaces
 \begin{equation*}
 B\St^{p}(\R, \espX):= \left\{f\in L_{\loc}^{p}(\R, \espX); \sup_{\xi\in \R}\left(\frac{1}{L}\int_{\xi}^{\xi+L}\norm{f(t)}^{p}dt\right)^{\frac{1}{p}}<\infty \right\} \forall L>0,
 \end{equation*}
 and
 \begin{equation*}
 B\We^{p}(\R, \espX):= \left\{f\in L_{\loc}^{p}(\R, \espX); \lim_{l\rightarrow\infty}\sup_{\xi\in \R}\left(\frac{1}{l}\int_{\xi}^{\xi+l}\norm{f(t)}^{p}dt\right)^{\frac{1}{p}}<\infty \right\},
 \end{equation*}
 coincide, see \cite{Andres_Bersani_Grande_06_Hierarchy}.
\end{enumerate}}
\subsection {Weyl almost periodic functions depending on a parameter }{\em
 \begin{enumerate}
 \item  We say that a parametric function $f: \mathbb{R}\times \espX\rightarrow \espX$ is {\em Weyl almost periodic} if, for every $ x \in \espX$, the function $f(.,x)$ is Weyl almost periodic, we denote by $\WeAP{p}(\R\times \espX, \espX)$ the space of such functions.

\item  A function $f: \mathbb{R}\times \espX\rightarrow \espX, (t,u)\mapsto f(t,u)$ with $f(.,u)\in L^{p}_{\loc}(\R, \espX)$ for each $u\in\espX$ is said to be {\em Weyl almost periodic in $t\in \mathbb{R}$ uniformly with respect compact subsets of $\espX$} if, for each $\epsilon>0$, there exists $l=l(\epsilon)>0$ and for all compacts $K$ in $\espX$, the set
\begin{equation*}
T_{\St^{p}_{l}}(\epsilon, f,K)=\left\{\tau\in\R;\sup_{u\in K}\sup_{\xi\in \mathbb{R}}\CCO{\frac{1}{l}\int_{\xi}^{\xi+l}\norm{f(t + \tau, u)-  f(t, u)}^{p}dt}^{\frac{1}{p}}< \epsilon\right\}
\end{equation*}
 is relatively dense. We denote by  $ \WeAP{p}_{K}(\R\times \espX, \espX)$ the set of such functions.
\end{enumerate}
 It is clear that $\WeAP{p}_{K}(\R\times \espX, \espX) \subset \WeAP{p}(\R\times \espX, \espX)$.
 }
\medskip

 {\em Let us introduce some notations due to Danilov \cite{danilov2006weyl}. We denote by $M_{p}^{*}(\R, \espX)$ the set of functions $f\in B\St^{p}(\R, \espX)$ such that
\begin{equation*}
\lim_{\delta\rightarrow0^{+}}\lim_{l_{0}\rightarrow\infty}\sup_{l>l_{0}}\sup_{\xi\in\R}\CCO{\frac{1}{l}\sup_{T\subseteq[\xi, \xi+l]:\,\, \varkappa(T)\leq\delta l }\int_{T}\norm{f(t)}^{p}dt}^{\frac{1}{p}} =0
\end{equation*}
holds, where $\varkappa$ is the Lebesgue measure on $\R$.

On the space $\espX$ we also consider the truncated norm $\norm{.}_{*} = \min\{1, \norm{.}\}$. Let $\WeAP{0}(\R, \espX) := \WeAP{1}(\R, (\espX, \norm{.}_{*}))$ denote the space of Weyl almost periodic functions $f : \R\rightarrow \espX$ (of degree $1$), when $\espX$ is endowed with  $\norm{.}_{*}$.
From Danilov \cite{danilov2006weyl}, we have
\begin{theo}[Danilov \cite{danilov2006weyl}]\label{prop:Danilov}
\begin{enumerate}
\item The following characterization holds true:
\begin{equation}\label{eq:charact_Danilov}
    \WeAP{p}(\R, \espX) = \WeAP{0}(\R, \espX)\cap M_{p}^{*}(\R, \espX).
\end{equation}
\item 
 If $f\in  \WeAP{p}(\R, \espX)$,
then  for  every $\epsilon, \delta > 0$ there exist a number $l>0$ and a compact set $K_{\epsilon, \delta} \subset \espX$ such that
\begin{equation*}
\sup_{\xi\in\R}\frac{1}{l}\varkappa\left\{t\in[\xi,\xi+l]: d\CCO{f(t), K_{\epsilon, \delta}}\geq\delta\right\}<\epsilon.
\end{equation*}
\end{enumerate}
\end{theo}

Note that, as for almost periodicity in Stepanov sense, which is a metric property  \cite[Remark~2.4]{bedouhene2017almost}, it can be easily shown that almost periodicity in Weyl sense is a metric property too. In fact, using the arguments in \cite[Remark~2.4]{bedouhene2017almost}, let $g=\exp\CCO{\sum_{n=2}^{+\infty}g_n}$, where $g_n$ is a $4n$-periodic function given by
\begin{equation*}
    g_n(t)=\beta_n\CCO{1-\frac{2}{\alpha_n}\abs{t-n}}
  \un{[n-\frac{2}{\alpha_n},n+\frac{2}{\alpha_n}]}(t),\quad
  t\in[-2n,2n],
\end{equation*}
where $\alpha_{n}\in ]0, \frac{1}{2}[$, and $\beta_{n}>0$. As shown in \cite{andres2012nonexistence}, for $\alpha_n=1/n^5$ and $\beta_n=n^3$, we have $\norm{g}_{\St^{1}}\geq \frac{\alpha_{n}\beta_{n}^{2}}{6}\rightarrow\infty$. Which means that $g$ is not $\St^{1}$-bounded. Hence $g$ cannot be in $\WeAP{1}(\R)$, . However, the function $g$ is almost periodic in Lebesgue measure, that is, it belongs to $\StAP{0}(\R):=\StAP{1}\big(\R,(\R, \norm{.}_{*})\big)$. In view of the inclusion $\StAP{0}(\R)\subset \WeAP{0}(\R)$, we deduce that $g$ belongs to $\WeAP{0}(\R)$. This shows also that the inclusion $\WeAP{p}(\R)\subset \WeAP{0}(\R)$ is strict, and that $M_{p}^{*}(\R, \espX)$ depends on the norm of the space $\espX$.

\section{Main results}
\subsection{A superposition theorem in $ \WeAP{p}(\R, \espX)$ }
  In this section, we establish a composition theorem for Weyl almost periodic functions. Our result is based on Theorem~\ref{prop:Danilov}. In the following, we assume that $\frac{1}{p}=\frac{1}{q}+ \frac{1}{r}$  with $p$, $q$ and $r \geq 1$.

 \theo \label{theorem:theorm superpostion} Let $f \in  \WeAP{p}_{K}(\R \times\espX, \espX)$. Assume that there exists a positive
function $L(.)\in B\We^{r}(\R)$, such that
\begin{equation}\label{For:Lip-condition}
 \norm{f(t, u)-f(t,v)} \leq L(t)\norm{u-v}\quad \forall t\in \R,\quad u, v \in \espX.
\end{equation}
Then, for every $x(.)\in \WeAP{q}(\R, \espX)$, we have
\begin{equation*}
f(., x(.)) \in \WeAP{p}(\R, \espX).
\end{equation*}
\proof{\em For a measurable set $A\subset\R$ and $l>0$, let us introduce the notation
\begin{equation*}
\overline{\varkappa}_{l}(A)= \sup_{\xi\in\R}\frac{1}{l}\varkappa\left\{t\in[\xi, \xi+l]\cap A\right\},
\end{equation*}
and let us denote by $A^{c}$ the complementary set of $A$. Fix  $\epsilon>0$. Let $x(.)\in \WeAP{q}(\R, \espX)$. In view of \eqref{eq:charact_Danilov}, we can find  $\delta_{\epsilon}>0$ and $l_{0}=l_{0}(\epsilon)>0$ such that, for every measurable set $A$, if $\overline{\varkappa}_{l_{0}}(A)\leq\delta_\epsilon$, then
\begin{equation}\label{for:unifo.inte}
\sup_{x\in\R}\CCO{\frac{1}{l_{0}}\int_{A\cap[x,x+l_{0}]}\norm{x(s)}^{p}ds}^\frac{1}{p}\leq\epsilon.
\end{equation}
Using the second item of Theorem~\ref{prop:Danilov}, for the above $\epsilon>0$ and the corresponding $\delta_\epsilon>0$, there exist a compact subset $K_{\epsilon}\subset\espX$ and a positive number $l_{1}=l_{1}(\epsilon)>0$ such that
\begin{equation}\label{for: condcorol}
\overline{\varkappa}_{l_{1}}(T_{\epsilon}) < \delta_\epsilon,
\end{equation}
where
\begin{equation*}
T_{\epsilon} = \left\{t\in\R: d(x(t), K_{\epsilon}) \geq \frac{\epsilon}{24\norm{L}_{\We^{r}}} \right\}.
\end{equation*}
Hence, we obtain by taking $l_{2}=\max(l_{0},l_{1})$,
\begin{equation}\label{for: unifo.integra}
\sup_{x\in\R}\CCO{\frac{1}{l_{2}}\int_{
T_{\epsilon}\cap[x,x+l_{2}]}\norm{x(s)}^{p}ds}^\frac{1}{p}\leq\frac{\epsilon}{24\norm{L}_{\We^{r}}}.
\end{equation}
The compactness of $K_{\epsilon}$ yields that for every $t\in\R$, there exists $x^{*}_{t}\in K_{\epsilon}$ such that
\begin{equation*}
d(x(t), K_{\epsilon}) = \norm{x(t)- x^{*}_{t}}.
\end{equation*}
Now, since $x(.) \in M_{p}^{*}(\R, \espX)$ and $\accol{x^{*}_{t}; \; t\in \R}\subset K_\epsilon$, we obtain by \eqref{for: unifo.integra} that,
\begin{equation}\label{for: equation sur T}
\sup_{\xi\in\R}\CCO{\frac{1}{l_{2}}\int_{T_{\epsilon}\cap[\xi,\xi+l_2]}\norm{x(t)- x^{*}_{t}}^{q}dt}^{\frac{1}{q}}<\frac{\epsilon}{24\norm{L}_{\We^{r}}}.
\end{equation}
We also have
\begin{equation}\label{for: equation sur T*}
\sup_{\xi\in\R}\CCO{\frac{1}{l_{2}}\int_{T_{\epsilon}^{c}\cap[\xi,\xi+l_{2}]}\norm{x(t) - x^{*}_{t}}^{q}dt}^{\frac{1}{q}}< \frac{\epsilon}{24\norm{L}_{\We^{r}}}.
\end{equation}
Using again the compactness of  $K_{\epsilon}$, we can find a finite sequence  $y_{1}, y_{2}, \cdots, y_{n} \in K_{\epsilon}$ such that
\begin{equation}\label{for: le compact K}
K_{\epsilon} \subset \bigcup_{i=1}^{n}B\CCO{y_{i}, \frac{\epsilon}{24\norm{L}_{\We^{r}}}}.
\end{equation}
Since $x^{*}_{t}\in K_{\epsilon}$, it follows that,  for every $t\in\R$, there exists $i(t)\in \left\{1, 2, \cdots, n\right\}$ such that
\begin{equation}\label{for: le compact k..}
\sup_{\xi\in\R}\CCO{\frac{1}{l_{1}}\int_{\xi}^{\xi+l}\norm{x^{*}_{t}-y_{i(t)}}^{q}dt}^{\frac{1}{q}}<\frac{\epsilon}{24\norm{L}_{\We^{r}}}.
\end{equation}
On the other hand, as $f \in \WeAP{p}_{K}(\R \times\espX, \espX)$ and $x(.)\in \WeAP{q}(\R, \espX)$, for the above  $\epsilon > 0$, we can chose $l_3=l_3(\epsilon)>0$ and a common relatively dense set $T_{\St^{p}_{l}}(\epsilon,x,f) \subset \R$ such that
\begin{equation}\label{For:equi-Weyl de x}
\sup_{\xi\in\R}\CCO{\frac{1}{l_3}\int_{\xi}^{\xi+l_3}\norm{x(t+\tau) - x(t)}^{q}dt}^\frac{1}{q}< \frac{\epsilon}{2\norm{L}_{\We^{r}}},
\end{equation}
and
\begin{equation}\label{For:equi-Weyl de f}
\sup_{\xi\in\R}\CCO{\frac{1}{l_3}\int_{\xi}^{\xi+l_3}\norm{f(t+\tau, v) - f(t, v)}^{p}dt}^\frac{1}{p}< \frac{\epsilon}{12n},
\end{equation}
for all  $\tau\in T_{\St^{p}_{l}}(\epsilon,x,f)$ and $v\in K_{\epsilon}$. Furthermore, by triangular inequality, we have, for $l=\max(l_2,l_3)$ and for every $\tau\in T_{\St^{p}_{l}}(\epsilon,x,f)$,
\begin{align*}
\sup_{\xi\in\R}\Big(\frac{1}{l}\int_{\xi}^{\xi+l}&\norm{f(t+\tau, x(t+\tau)) - f(t,
x(t))}^{p}dt\Big)^{\frac{1}{p}}\\ \leq&
\sup_{\xi\in\R}\CCO{\frac{1}{l}\int_{\xi}^{\xi+l}\norm{f(t+\tau, x(t+\tau))
- f(t+\tau,x(t))}^{p}dt}^{\frac{1}{p}}\\&+ \sup_{\xi\in\R}\CCO{\frac{1}{l}\int_{\xi}^{\xi+l}\norm{f(t+\tau, x(t)) -
f(t, x(t))}^{p}dt}^{\frac{1}{p}}\\ =:& I_{1} + I_{2}.
\end{align*}
Let us estimate $I_{1}$ and $I_{2}$. For $I_{1}$, by \eqref{For:Lip-condition}, \eqref{For:equi-Weyl de x} and H\"older's inequality ($\frac{r}{p},\frac{q}{p}$), we have, for every  $\tau\in T_{\St^{p}_{l}}(\epsilon,x,f)$,
\begin{align*}
    I_{1} = \sup_{\xi\in\R}&\CCO{\frac{1}{l}\int_{\xi}^{\xi+l}\norm{f(t+\tau, x(t+\tau))
- f(t+\tau,x(t))}^{p}dt}^{\frac{1}{p}}\\ &\leq \norm{L}_{\We^{r}} \sup_{\xi\in\R}\CCO{\frac{1}{l}\int_{\xi}^{\xi+l}\norm{ x(t+\tau) -
x(t)}^{q}dt}^\frac{1}{q}< \frac{\epsilon}{2}.
\end{align*}
For $I_{2}$, for all $\tau\in T_{\St^{p}_{l}}(\epsilon,x,f)$, we have
\begin{align*}
I_{2}  \leq& \sup_{\xi\in\R}\CCO{\frac{1}{l}\int_{T_{\epsilon}\cap[\xi,\xi+l]}\norm{f(t+\tau, x(t)) - f(t, x(t))}^{p}dt}^{\frac{1}{p}}
\\&+ \sup_{\xi\in\R}\CCO{\frac{1}{l}\int_{T_{\epsilon}^{c}\cap[\xi,\xi+l]}\norm{f(t+\tau, x(t)) - f(t, x(t))}^{p}dt}^{\frac{1}{p}} :=\,\,\,\, I_{2}^{1} +
I_{2}^{2}.
\end{align*}
First, we consider $I_{2}^{1}$. Using triangular inequality, the condition \eqref{For:Lip-condition}, \eqref{for: equation sur T} and H\"older's inequality ($\frac{p}{r}, \frac{p}{q}$), we get, for all $\tau\in T_{\St^{p}_{l}}(\epsilon,x,f)$,
\begin{align*}
 I_{2}^{1} &= \sup_{\xi\in\R}\CCO{\frac{1}{l}\int_{T_{\epsilon}\cap[\xi,\xi+l]}\norm{f(t+\tau, x(t)) - f(t, x(t))}^{p}dt}^{\frac{1}{p}}\\ \leq&
 \frac{\epsilon}{24} + \sup_{\xi\in\R}\CCO{\frac{1}{l}\int_{\xi}^{\xi+l} \norm{f(t+\tau, x^{*}_t) - f(t, x^{*}_t)}^{p}dt}^{\frac{1}{p}} +\frac{\epsilon}{24}\\:=& \frac{\epsilon}{12}+ I^{*}.
\end{align*}
By \eqref{For:equi-Weyl de f}, for each $\tau\in T_{\St^{p}_{l}}(\epsilon,x,f)$ and $i = 1, 2, \cdots, n$, we obtain that
\begin{equation}\label{for: double8}
\sup_{\xi\in\R}\CCO{\frac{1}{l}\int_{\xi}^{\xi+l}\norm{f(t+\tau, y_{i}) - f(t, y_{i})}^{p}dt}^\frac{1}{p}< \frac{\epsilon}{12n}.
\end{equation}
Using \eqref{For:Lip-condition} and \eqref{for: le compact k..}, it follows that, for every $\tau\in T_{\St^{p}_{l}}(\epsilon,x,f)$,
\begin{align*}
I^{*} \leq& 2\norm{L}_{\We^{r}}\sup_{\xi\in\R}\CCO{\frac{1}{l}\int_{\xi}^{\xi+l}\norm{y_{i(t)} - x^{*}_{t}}^{q}dt}^{\frac{1}{q}} \\&+
\sup_{\xi\in\R}\CCO{\frac{1}{l}\int_{\xi}^{\xi+l}\norm{f(t+\tau, y_{i(t)}) - f(t, y_{i(t)})}^{p}dt}^{\frac{1}{p}}\\ \leq&
 \frac{\epsilon}{12}+  \sum_{j=1}^{n_{l}} \sup_{\xi\in\R}\CCO{\frac{1}{l}\int_{\xi}^{\xi+l}\norm{f(t+\tau, y_{j}) - f(t, y_{j})}^{p}dt}^\frac{1}{p}< \frac{\epsilon}{6}.
\end{align*}
We thus have, for all $\tau\in T_{\St^{p}_{l}}(\epsilon,x,f)$,
\begin{equation*}
I_{2}^{1} = \sup_{\xi\in\R}\CCO{\frac{1}{l}\int_{T_{\epsilon}\cap[\xi,\xi+l]}\norm{f(t+\tau, x(t)) - f(t, x(t))}^{p}dt}^{\frac{1}{p}} < \frac{\epsilon}{4}.
\end{equation*}
For $I_{2}^{2}$, by \eqref{For:Lip-condition}, \eqref{for: equation sur T*}, \eqref{for: le compact k..}, \eqref{for: double8}  and H\"older's inequality, we follow the same steps than for $I_{2}^{1}$, we obtain, for every $\tau\in T_{\St^{p}_{l}}(\epsilon,x,f)$,
\begin{equation*}
 I_{2}^{2} = \sup_{\xi\in\R}\CCO{\frac{1}{l}\int_{T_{\epsilon}^{c}\cap[\xi,\xi+l]}\norm{f(t+\tau, x(t)) - f(t, x(t))}^{p}dt}^{\frac{1}{p}}< \frac{\epsilon}{4}.
\end{equation*}
Then, for all $\tau\in T_{\St^{p}_{l}}(\epsilon,x,f)$, we get that
\begin{equation*}
I_{2} = \sup_{\xi\in\R}\CCO{\frac{1}{l}\int_{\xi}^{\xi+l}\norm{f(t+\tau, x(t)) - f(t, x(t))}^{p}dt}^{\frac{1}{p}} < \frac{\epsilon}{2}.
\end{equation*}
Finally, combining $I_{1}$ and $I_{2}$, we obtain, for every $\tau\in T_{\St^{p}_{l}}(\epsilon,x,f)$,
\begin{equation*}
 \sup_{\xi\in\R}\CCO{\frac{1}{l}\int_{\xi}^{\xi+l}\norm{f(t+\tau, x(t+\tau)) - f(t, x(t))}^{p}dt}^{\frac{1}{p}}< \epsilon,
\end{equation*}
that is,  $f(., x(.))\in \WeAP{p}(\R, \espX)$.} \finpr \medskip

{\em Now, we discuss the Weyl almost periodic solutions to linear differential equations.}
 \subsection{Existence and uniqueness of  bounded Weyl almost periodic solution to abstract linear differential equation}
{\em We investigate now the Weyl almost periodic solutions to
 \begin{equation}\label{For: equation lineaire}
 u^{\prime }\left( t\right) =Au(t) +f(t),
\end{equation}
where $A:D\left( A\right) \subset \espX\rightarrow \espX$ is a linear operator  which generates an exponentially stable  $C_{0}$-semigroup $(T(t))_{t\geq 0}$   on a Banach space $\espX$, and $f :\mathbb{R} \rightarrow \espX $ is Weyl almost periodic.

We recall that a function $u:\R\rightarrow\espX$ is a {\em mild solution} to \eqref{For: equation lineaire} if
\begin{equation*}
u(t) = T(t-a)u(a) + \int_{a}^{b}T(t-s)f(s)ds,\quad t\geq s.
\end{equation*}
Let us also recall that a semigroup $(T(t))_{t\in\R}$ of linear operators is {\em exponentially stable} if there exist $\delta>0$ and $M\geq 1$ such that
\begin{equation}\label{For: hyp.exponentially stable}
\norm{T(t)}\leq M e^{-\delta t}.
\end{equation}

The main motivation of this study is that, contrarily the Stepanov case (see \cite{andres2012nonexistence} and \cite{bedouhene2017almost}), where it is shown the  nonexistence of purely Stepanov almost periodic solutions for ordinary differential equations with Stepanov almost periodic coefficients, we can find, in the Weyl case, a bounded and purely Weyl almost periodic mild solution to Equation \eqref{For: equation lineaire}. The problem is closely related to Bohl-Bohr theorem which fails in the Weyl case. Indeed, the following example shows that there exist Weyl almost periodic functions such that their integrals are bounded but not Bohr almost periodic, nor Stepanov almost periodic. However, they are (purely) Weyl almost periodic.}} 
 {\em \begin{example}
    {\em We consider the following function $f:\R\rightarrow\R$ defined by
  \begin{equation}\label{eq:fonction_counter_example}
f(t) = \left\{%
\begin{array}{ll}
    1  & \hbox{for} \quad 0<t<\frac{1}{2}; \\
    0  & \hbox{elsewhere},
\end{array}%
\right.
\end{equation}
which is purely Weyl almost periodic (see \cite[p.145]{Andres_Bersani_Grande_06_Hierarchy}). Let $F:\R\rightarrow\R$ with
  \begin{equation*}
F(t) =\int_{-\infty}^{t}f(s)ds= \left\{%
\begin{array}{lll}
    0  & \hbox{for} \quad t\leq 0; \\
    t  & \hbox{for} \quad 0<t<\frac{1}{2}; \\
    \frac{1}{2}  & \hbox{elsewhere},
\end{array}%
\right.
\end{equation*}
be a primitive of $f$. Clearly $F$ is bounded on $\R$. Moreover, for every $t_{1}, t_{2}\in\R$, we have
\begin{equation*}
\abs{F(t_{1})-F(t_{2})}\leq \abs{t_{1}-t_{2}}.
\end{equation*}
Thus $F$  is uniformly continuous. So by \cite[Proposition~4]{Radova_04}, $F$ is Weyl almost periodic.

 On the other hand, for $\epsilon < \frac{1}{4}$, there exists $t\in[0,\frac{1}{4}]$ such that, for every $\tau > \epsilon$, we have
\begin{equation*}
|F(t+\tau)-F(t)|\geq \frac{1}{4}>\epsilon.
\end{equation*}
So, the function $F$ is not Bohr almost periodic. We can show more, that is $F$ is not Stepanov almost periodic. Indeed, if $F$ is Stepanov almost periodic, thanks to its uniform continuity, we get would have that $F$ was Bohr almost periodic (Bochner's Theorem, see. e.g \cite{Andres_Bersani_Grande_06_Hierarchy}, \cite[p.82]{besicovitoh1932almost}, \cite[Lemme.4 p.34]{levitan1982almost} and \cite[Th 6.16]{Corduneanu}). A contradiction. Hence $F$ is not Stepanov almost periodic.}
 \end{example}}
 {\em \begin{example}
  {\em Let the scalar differential equation $x'=-x+f(t)$, $t\in\R$, with $f$ given by \eqref{eq:fonction_counter_example}. Its unique bounded solution is given by
   \begin{equation}
     x(t)=\left\{
     \begin{array}{lll}
    0  & \hbox{for} \quad t\leq 0; \\
    1-\exp(-t)  & \hbox{for} \quad 0<t<\frac{1}{2}; \\
    (\sqrt{\exp(1)}-1)\exp(-t)  & \hbox{elsewhere}.
\end{array}%
          \right.
   \end{equation}
   The solution $x$ is clearly uniformly continuous. Moreover, by
   \cite[Proposition~4]{Radova_04}, $x$ is Weyl almost periodic. Now,
   we can get more, namely that the solution $x$ is purely Weyl almost periodic. Indeed, for every $\epsilon>0$, the following inequality $$(\sqrt{\exp(1)}-1)\exp(-t) < \epsilon$$ holds for all $t\geq -\ln\CCO{\frac{\epsilon}{\sqrt{\exp(1)}-1}}$. Now, let us choose $\epsilon=\frac{\sqrt{\exp(1)}}{2\CCO{\sqrt{\exp(1)}-1}}$. Then, every $\tau$ such that $\tau+\frac{1}{2}>-\ln\CCO{\frac{\epsilon}{\sqrt{\exp(1)}-1}}$, satisfies
   \begin{equation*}
    x\CCO{\frac{1}{2}}-x\CCO{\tau+\frac{1}{2}}=2\epsilon-(\sqrt{\exp(1)}-1)\exp(-\tau-\frac{1}{2}) \geq 2\epsilon-\epsilon=\epsilon.
   \end{equation*}
   Which means that $x$ is not almost periodic, nor Stepanov almost periodic.}
 \end{example}}
{\em \theo \label{Theorem: linear} There exists a unique bounded mild solution to \eqref{For: equation lineaire}. This solution is Weyl almost periodic of degree $p$ and is given by
\begin{equation}\label{For: la solution u}
u(t) = \int_{-\infty}^{t}T(t-s)f(s)ds.
\end{equation}
\proof
 {\em First, we show that the function $u(t) = \int_{-\infty}^{t}T(t-s)f(s)ds$  is well defined, i.e, $\lim_{a\rightarrow-\infty}\int_{a}^{t}T(t-s)f(s)ds$
 exists. We set
\begin{equation*}
    G(a) = \int_{a}^{t}T(t-s)f(s)ds.
\end{equation*}
Let us show that $(G(a))_{a\in\R}$ is Cauchy at $-\infty$. If $a< b$, then by \eqref{For: hyp.exponentially stable}, we have
\begin{equation*}
\norm{G(b)-G(a)} \leq \int_{a}^{b}\norm{T(t-s)}\norm{f(s)}ds \leq M\int_{a}^{b}e^{-\delta(t-s)}\norm{f(s)}ds.
\end{equation*}
Since $f$ is Weyl almost periodic, we have $f\in B\St^{p}_{L}(\R, \espX)$, for every $L>0$.
Thus, there exists $N\in\N$, such that, by H\"older's inequality ($\frac{1}{p}, \frac{1}{q}$), we have
\begin{align*}
\norm{G(b)-G(a)}& \leq M\sum_{k=0}^{N}\int_{b- (k+1)L}^{b- kL}e^{-\delta(t-s)}\norm{f(s)}ds\\& \leq M\sum_{k=0}^{N}\CCO{\int_{b- (k+1)L}^{b- kL}e^{-\delta q(t-s)}ds}^{\frac{1}{q}}\CCO{\int_{b- (k+1)L}^{b- kL}\norm{f(s)}^{p}ds}^{\frac{1}{p}}\\&\leq ML\norm{f}_{\St_{L}^{p}}e^{-\delta q(t-b)}\CCO{\frac{1-e^{-\delta q(N+1)L}}{1-e^{-\delta qL}}}.
\end{align*}
Consequently,
\begin{equation*}
\lim_{b \rightarrow-\infty} \norm{G(b)-G(a)}\leq \lim_{b \rightarrow-\infty} ML\norm{f}_{\St_{L}^{p}} e^{-\delta q(t-b)}\CCO{\frac{1-e^{-\delta q(N+1)L}}{1-e^{-\delta qL}}}=0.
\end{equation*}
Thus the limit
\begin{equation*}
 \lim_{a\rightarrow-\infty}\int_{a}^{t}T(t-s)f(s)ds
\end{equation*}
exists. This means that $u$ is well-defined. Moreover, it is easy to see that $u$ given by \eqref{For: la solution u} is a mild solution to \eqref{For: equation lineaire}.

In the following, we show that $u$ is bounded. By H\"older's inequality and \eqref{For: hyp.exponentially stable}, we obtain
\begin{equation*}
\norm{u(t)}\leq \int_{-\infty}^{t}\norm{T(t-s)}\norm{f(s)}ds \leq M\sum_{k=0}^{\infty} e^{-\delta k}\CCO{\int_{t-k-1}^{t-k}\norm{f(s)}^{p}ds}^{\frac{1}{p}}.
\end{equation*}
Then, we have
\begin{equation*}
\sup_{t\in \R}\norm{u(t)} \leq  M\sum_{k=0}^{\infty}e^{-\delta k}\sup_{t\in \R}\CCO{\int_{t-k-1}^{t-k}\norm{f(s)}^{p}ds}^{\frac{1}{p}}.
\end{equation*}
Since Weyl almost periodic functions of degree $p$ are $\St^{p}$-bounded, we have
\begin{equation*}
\sup_{t\in\R}\norm{u(t)}_{\infty}\leq M\norm{f}_{\St^{p}}\sum_{k=0}^{\infty} e^{-\delta k}<\infty.
\end{equation*}
Thus $u$ is bounded.

Now, we prove that the solution to \eqref{For: equation lineaire} is Weyl almost periodic. Using a change of variable  and H\"older's inequality $(\frac{1}{p}+\frac{1}{q}=1)$, we have any $\tau\in\R$,
\begin{align*}
\norm{u(t+\tau)-u(t)}&\leq M\int_{-\infty}^{0}e^{\delta s}\norm{f(t+s+\tau)-f(t+s)} ds  \\ &\leq M\CCO{\frac{2}{\delta
q}}^{\frac{1}{q}}\CCO{\int_{-\infty}^{0}e^{\frac{\delta ps}{2}}\norm{f(t+s+\tau)-f(t+s)}^{p} ds}^{\frac{1}{p}} .
\end{align*}
Let $l>0$. Using Fubini's theorem, we get, for every $\tau\in\R$,
\begin{align*}
\sup_{\xi\in \R}\ &\CCO{\frac{1}{l}\int_{\xi}^{\xi+l}\norm{u(t+\tau)-u(t)}^{p}dt}^{\frac{1}{p}}\\&\leq M\CCO{\frac{2}{\delta q}}^{\frac{1}{q}}
\sup_{\xi\in \R}\ \CCO{\frac{1}{l}\int_{\xi}^{\xi+l}  \CCO{\CCO{\int_{-\infty}^{0}e^{\frac{\delta ps}{2}}\norm{f(t+s+\tau)-f(t+s)}^{p} ds}^{\frac{1}{p}}}^{p}dt}^{\frac{1}{p}}\\&\leq M\CCO{\frac{2}{\delta q}}^{\frac{1}{q}} \CCO{\int_{-\infty}^{0}e^{\frac{\delta ps}{2}}\sup_{\xi\in \R}\frac{1}{l}\int_{\xi}^{\xi+l}\norm{f(t+s+\tau)-f(t+s)}^{p} dt\,ds}^{\frac{1}{p}}.
\end{align*}
Since $f$ is Weyl almost periodic, for any $\epsilon > 0$, we can choose positive number $l=l(\epsilon)$ and there exists a relatively dense set $T_{\St^{p}_{l}}(\epsilon,f)$ as in \eqref{for: definition equi-Weyl}, such that
\begin{equation*}
\CCO{\int_{-\infty}^{0}e^{\frac{\delta ps}{2}}\sup_{\xi\in \R}\frac{1}{l}\int_{\xi}^{\xi+l}\norm{f(t+s+\tau)-f(t+s)}^{p} dt\,ds}^{\frac{1}{p}}\leq \CCO{\frac{2}{\delta p}}^{\frac{1}{p}}\epsilon\quad\forall\tau\in T_{\St^{p}_{l}}(\epsilon,f).
\end{equation*}
Hence,
\begin{equation*}
\sup_{\xi\in \R}\ \CCO{\frac{1}{l}\int_{\xi}^{\xi+l}\norm{u(t+\tau)-u(t)}^{p}dt}^{\frac{1}{p}}\leq C\epsilon,
\end{equation*}
where $C= M\CCO{\frac{2}{\delta q}}^{\frac{1}{q}}\CCO{\frac{2}{\delta p}}^{\frac{1}{p}}$. Summing up, we have shown that, for every $\epsilon>0$, there corresponds a positive number $l=l(\epsilon)$ and there exists a relatively dense set $T_{\St^{p}_{l}}(\epsilon,u)=T_{\St^{p}_{l}}(\frac{\epsilon}{C},f)$ such that
\begin{equation*}
 \sup_{\xi\in \R}\
 \CCO{\frac{1}{l}\int_{\xi}^{\xi+l}\norm{u(t+\tau)-u(t)}^{p}dt}^{\frac{1}{p}} < \epsilon;  \quad\quad \forall \tau \in T_{\St^{p}_{l}}(\epsilon,u).
\end{equation*}
Which implies the Weyl-almost-periodicity of the solution $u$.

Uniqueness of the solution follows from the same arguments as those in \cite{Ding_Long_Nguerekata2011,long2011composition}. We omit the details. \finpr}
\subsection{Existence and uniqueness of  bounded Weyl almost periodic solution to abstract semilinear differential equation }
 {\em In this section we study the existence and uniqueness of the  Weyl almost periodic mild solution to equation
 \begin{equation}\label{for: equation principal}
 u'(t) = Au(t) + f(t, u(t)),
 \end{equation}
 where  $A:D\left( A\right) \subset \espX\rightarrow \espX$ is unbounded linear, of domain $D(A)$ dense in $\espX$, and $f:\R\times\espX\rightarrow\espX$ measurable function.

We consider the following hypothesis:
\begin{itemize}
\item[(H1)] The operator $A$ generates a $C_{0}-$semigroup $(T(t))_{t\geq0}$, exponentially stable.
\item[(H2)] $f\in\WeAP{p}_{K}(\R\times \espX,\espX)$.
\item[(H3)]The function $f$ is L(.)-Lipschitz i.e., there exists a  nonnegative function $L(.)
\in B\We^{p}(\R)$ such that
\begin{equation*}
\norm{f(t, u)-f(t,v)} \leq L(t)\norm{u-v};\quad \forall t\in \R,\quad \forall\, u,
v \in \espX.
\end{equation*}
\end{itemize}
}
{\em \theo \label{theorem: principale} Let the assumptions {\em(H1)-(H3)} be fulfilled. For $p \geq 2$,
the equation \eqref{eq:EDO semi-lineaire} has a  unique  mild solution in $\CB{(\R, \espX)}$ given by
\begin{equation}\label{For: la solution de l'equation semi-linear}
u(t) = \int_{-\infty}^{t}T(t-s)f(s,u(s))ds,
\end{equation}
provided that
\begin{equation}\label{For: L{Sp p>=2}}
\norm{L}_{\St^{p}}< \CCO{\frac{\delta q}{M^{q}}}^{\frac{1}{q}}\CCO{\frac{e^{\delta}-1}{e^{\delta}}},
\end{equation}
where $\frac{1}{p}+ \frac{1}{q} = 1$.

If, furthermore
\begin{equation}\label{For: L{Sp p>2}}
\norm{L}_{\St^{p}}<\frac{p\delta}{M^{p}2^{p+1}} \CCO{\frac{e^{\frac{p\delta}{4}} - 1}{e^{\frac{p\delta}{4}}}}\CCO{\frac{p\delta}{2p-4}}^{p-2}\hbox{for}\,\,\,p>2,
\end{equation}
or
\begin{equation}\label{For: L{Sp p=2}}
\norm{L}_{\St^{2}}<\frac{p\delta}{8M^{2}}\CCO{\frac{e^{\frac{2\delta}{4}} - 1}{e^{\frac{2\delta}{4}}}}\quad\quad\quad\quad\quad\quad\hbox{for}\,\,\,\, p=2,
\end{equation}
then, $u$ is Weyl almost periodic.\medskip\\
{\em Before giving the proof of Theorem \ref{theorem: principale}, we need the following technical results}.
{\em \lemma \label{Lemma: Kame} {\em(\cite{Kamenskii_Mikhail_and_Mellah_Omar_and_Raynaud_de_Fitte_15})} Let $g: \R\rightarrow \R$ be a
continuous function such that, for every $t\in\R$,
\begin{equation}\label{For: Lemma}
0\leq g(t)\leq \alpha(t)+
\beta_{1}\int_{-\infty}^{t}e^{-\delta_{1}(t-s)}g(s)ds+...+
\beta_{n}\int_{-\infty}^{t}e^{-\delta_{1}(t-s)}g(s)ds,
\end{equation}
for some locally integrable function $\alpha:\R\rightarrow\R$, and
for some constants $\beta_{1},...,\beta_{n}\geq0$, and some constants
$\delta_{1},...,\delta_{n}>\beta$, where $\beta =
\sum_{i=1}^{i=n}\beta_{i}$. We assume that the integrals in the right hand side of \eqref{For: Lemma} are convergent. Let $\delta = \min_{1\leq i\leq n}\delta_{i}$. Then, for every $\gamma\in ]0, \delta-\beta]$ such that $\int_{-\infty}^{0}e^{\gamma s}\alpha(s)ds$ converges, we have, for every $t\in\R$,
\begin{equation*}
g(t)\leq \alpha(t)+
\beta\int_{-\infty}^{t}e^{-\gamma(t-s)}\alpha(s)ds.
\end{equation*}
In particular, if $\alpha$ is constant, we have
\begin{equation*}
g(t)\leq \alpha\frac{\delta}{\delta-\beta}.
\end{equation*}
{\em \prop \label{proposition: principale}} Let $f$ be in $\WeAP{p}_{K}(\R\times\espX, \espX)$ $(2\leq p<\infty)$ and satisfying {\rm (H1)}. Let $u\in\CB{(\R, \espX)}$. Then, for every $\epsilon>0$ there exists $l=l(\epsilon)>0$ such that  for every $\tau\in T_{\St^{p}_{l}}(\epsilon,f,K)$, we have
\begin{equation}\label{For: proposition}
\int_{-\infty}^{0}e^{\delta_1 s} \frac{1}{l}\int_{s}^{s+l} \norm{f(t+\tau, u(t))- f(t, u(t))}^{p}dt\,ds<\epsilon^{p}
\end{equation}
and
\begin{equation}\label{For: proposition*}
 \int_{-\infty}^{0}e^{\gamma r}\int_{-\infty}^{r}e^{-\delta_1(r-s)} \frac{1}{l}\int_{s}^{s+l} \norm{f(t+\tau, u(t))- f(t, u(t))}^{p}dt\,ds\,dr<\epsilon^{p},
\end{equation}
for all $\delta_1>0$ and $\gamma>0$.

\begin{proof}
{\em Let fix $\delta_1>0$ and $\gamma>0$. For $s\in\R$, $\xi \in \R$, $l>0$ and $\tau\in\R$, let
\begin{eqnarray*}
h_{\tau, l}(s)&=& \frac{1}{l}\int_{s}^{s+l} \norm{f(t+\tau, u(t))- f(t, u(t))}^{p}dt,\\
\alpha_{\tau, l}(\xi)&=&\int_{-\infty}^{\xi}e^{-\delta_1(\xi-s)}h_{\tau, l}(s)ds.
\end{eqnarray*}
%
Inequalities \eqref{For: proposition} and \eqref{For: proposition*} become $\alpha_{\tau, l}(0)<\epsilon^{p}$ and $\int_{-\infty}^{0}e^{\gamma r}\alpha_{\tau, l}(r)dr<\epsilon^{p}$, respectively.

First, we deal with the estimation \eqref{For: proposition}. Using the boundedness of both $u$ and $f$ in the Weyl sense ($f$ is $\We^{p}$-bounded), we have for every $\tau\in\R$, $l>0$ and for every $s\in]-\infty,0]$,
\begin{align*}
h_{\tau, l}(s)=&\frac{1}{l}\int_{s}^{s+l} \norm{f(t+\tau, u(t))- f(t, u(t))}^{p}dt\\ \leq& 2^{2p-2}\CCO{\frac{1}{l}\int_{s}^{s+l} \norm{f(t+\tau, u(t))- f(t+\tau, 0)}^{p}dt+\frac{1}{l}\int_{s}^{s+l} \norm{f(t+\tau, 0)}^{p}dt} \\&+ 2^{2p-2}\CCO{\frac{1}{l}\int_{s}^{s+l} \norm{f(t, 0)}^{p}dt+\frac{1}{l}\int_{s}^{s+l} \norm{f(t, 0)- f(t, u(t))}^{p}dt} \\ \leq & 2^{2p-1}\CCO{\norm{L(.)}_{\We^{p}}^{p}\norm{u(.)}^{p}_\infty + \norm{f(.,0)}_{\We^{p}}^{p}} < +\infty.
\end{align*}
The function $s\mapsto e^{\delta_{1}s} h_{\tau, l}(s)$ is bounded and continuous on $]-\infty,0]$, hence
\begin{equation}\label{eq:estimate alpha propo}
\alpha_{\tau, l}(0)=\int_{-\infty}^{0}e^{\delta_{1}s} h_{\tau, l}(s)ds = \sum_{n=0}^{\infty} \int_{-n-1}^{-n}e^{\delta_{1}s} h_{\tau, l}(s)ds \leq \sum_{n=0}^{\infty} e^{-\delta_{1}n}\CCO{\max_{s\in[-n-1, -n]}h_{\tau, l}(s)}.
\end{equation}
Due to the uniform boundedness of the sequence $\CCO{\max_{s\in[-n-1, -n]} h_{\tau, l}(s)}_{n}$, the series in \eqref{eq:estimate alpha propo} is uniformly convergent. Thus, for every $\epsilon>0$, we can find an integer $N_{1}=N_{1}(\epsilon)\geq 1$ such that for every $\tau \in \R$, $l>0$, and $\delta_1>0$, we have
\begin{equation}\label{eq:estimate alpha propo_first_term}
\sum_{n> N_1} e^{-\delta_{1}n}\CCO{\max_{s\in[-n-1, -n]} h_{\tau, l}(s)} <\frac{\epsilon^{p}}{2}.
\end{equation}
Moreover, we claim that there exists $l=l(\epsilon)>0$ such that every $\tau \in T_{\St^{p}_{l}}(\epsilon,f,K)$ satisfies 
\begin{equation}\label{eq:estimate alpha propo_rest_term}
\sum_{n=0}^{N_1} e^{-\delta_{1}n}\CCO{\max_{s\in[-n-1, -n]}h_{\tau, l}(s)} < \frac{\epsilon^{p}}{2}.
\end{equation}
In fact, for the above  $\epsilon>0$, let $l=l(\epsilon)>0$ provided
by the condition that $T_{\St^{p}_{l}}(\epsilon,f,K)$ is relatively dense. In view of the continuity of the mapping  $s\mapsto h_{\tau, l}(s)$ on the compact set $[-n-1, -n]$, there exists  $s_{n}^{*}\in [-n-1, -n]$ such that, for every $\tau\in\R$ and $\delta_1>0$, $$\max_{s\in[-n-1, -n]}h_{\tau, l}(s)=h_{\tau, l}(s_{n}^{*}).$$
Put $s^{*}= s^{*}_{\epsilon}:=\displaystyle \argmax_{n\in\{0,\cdots, N_1\}}h_{\tau,l}(s_{n}^{*})$. We have then
\begin{align*}
 \sum_{n=0}^{N_1} e^{-\delta_{1}n}\max_{s\in[-n-1, -n]}h_{\tau, l}(s)\leq
    \sum_{n=0}^{N_1}e^{-\delta_{1}n} h_{\tau,l}(s^{*}).
\end{align*}
Now, since $u$ is continuous, we have that $K^*=u\left([s^{*}, s^{*}+ l]\right)$ is a compact set that depends only on  $\epsilon$. Hence, there exist $u_{1},\cdots, u_{k} \in K_{1}$ such that
\begin{equation*}
K^*\subset\bigcup_{i= 1}^{i= k}B\left(u_{i}, \frac{\epsilon}{2\sqrt[p]{6}\norm{L}_{\We^{p}}}\right).
\end{equation*}
It follows that, for each $t\in [s^{*}, s^{*}+l]$, there exists $i(t)\in \{1,2,\cdots,k\}$ such that
\begin{equation}\label{for:compact K1}
\norm{u(t)-u_{i(t)}}< \frac{\epsilon}{2\sqrt[p]{6}\norm{L}_{\We^{p}}}.
\end{equation}
On the other hand, since $f\in\WeAP{p}_{K}(\R\times \espX, \espX)$, every $\tau\in T_{\St^{p}_{l}}(\epsilon,f,K)$ is such that
\begin{equation}\label{for:alpha}
\frac{1}{l}\int_{s^{*}}^{s^{*}+l}\norm{f(t+\tau, u_{i(t)})-f(t, u_{i(t)})}^{p}dt < \frac{\epsilon^{p}}{2^{p+1}3k}.
\end{equation}
To combine the previous arguments, we use the triangle inequality:
\begin{align}
h_{\tau, l}(s^{*}) \leq&2^{p}\frac{1}{l}\int_{s^{*}}^{s^{*}+l} \norm{f(t+\tau, u(t))- f(t+\tau, u_{i(t)})}^{p}dt+ 2^{p}\frac{1}{l}\int_{s^{*}}^{s^{*}+l} \norm{f(t+\tau, u_{i(t)})- f(t, u_{i(t)})}^{p}dt \notag\\
&+ 2^{p}\frac{1}{l}\int_{s^{*}}^{s^{*}+l} \norm{f(t, u_{i(t)})- f(t, u(t))}^{p}dt \notag\\
=:& 2^{p}\CCO{A_{1}+A_{2}+A_{3}}.\label{eq:somme h_tau,l}
\end{align}
Let us estimate $A_{1}$. By (H2) and \eqref{for:compact K1}, we have
\begin{align}\label{eq:estimate A1 propo}
A_{1} &= \frac{1}{l}\int_{s^{*}}^{s^{*}+l} \norm{f(t+\tau, u(t))- f(t+\tau, u_{i(t)})}^{p}dt \notag\\
&\leq \frac{1}{l}\int_{s^{*}}^{s^{*}+l} \norm{L(t+\tau)}^{p} \norm{u(t)-
u_{i(t)}}^{p}dt\leq \CCO{\frac{\epsilon}{2\sqrt[p]{6}\norm{L}_{\We^{p}}}}^{p}\norm{L}_{\We^{p}}^{p}< \frac{\epsilon^{p}}{2^{p+1}3}.
\end{align}
We follow the same steps as for $A_{1}$, we obtain
\begin{equation}\label{eq:estimate A2 propo}
A_{3}\leq \CCO{\frac{\epsilon}{2\sqrt[p]{6}\norm{L}_{\We^{p}}}}^{p}\norm{L}_{\We^{p}}^{p}< \frac{\epsilon^{p}}{2^{p+1}3}.
\end{equation}
Now, we estimate $A_{2}$. By \eqref{for:alpha}, we have
\begin{align}\label{eq:estimate A3 propo}
A_{2}
\leq \sum_{j=1}^{k}\frac{1}{l}\int_{s^{*}}^{s^{*}+l} \norm{f(t+\tau, u_{j})- f(t, u_{j})}^{p}dt \leq \frac{\epsilon^{p}}{2^{p+1}3}.
\end{align}
Taking into account \eqref{eq:somme h_tau,l}-\eqref{eq:estimate A3 propo}, we obtain that
\begin{equation}\label{eq:estimation_final_h_tau,l}
    h_{\tau, l}(s^{*}) \leq \frac{\epsilon^{p}}{2},\, \text{ for every } \tau\in T_{\St^{p}_{l}}(\epsilon,f,K).
\end{equation}
This implies that
\begin{equation*}
\sum_{n=0}^{N_{1}} e^{-\delta_{1}n}\CCO{\max_{s\in[-n-1, -n]} h_{\tau,l}(s)} \leq 
(1-e^{-\delta_{1}})\frac{\epsilon^{p}}{2}\leq \frac{\epsilon^{p}}{2}.
\end{equation*}
Having disposed the previous preliminaries steps, we now return to estimation \eqref{eq:estimate alpha propo}, we get in view of \eqref{eq:estimate alpha propo_first_term} and \eqref{eq:estimate alpha propo_rest_term} that for every  $\epsilon>0$, there is $l=l(\epsilon)>0$ (provided by the condition $f\in \StAP[l]{p}_{K}(\R\times\espX, \espX)$) such that any $\tau \in T_{\St^{p}_{l}}(\epsilon,f,K)$ verifies
\begin{equation}\label{eq:estimate_alpha_0_last}
\alpha_{\tau, l}(0)=\int_{-\infty}^{0}e^{\frac{p\delta}{4}s} \frac{1}{l}\int_{s}^{s+l} \norm{f(t+\tau, u(t))- f(t, u(t))}^{p}dt\,ds \leq\epsilon^{p}.
\end{equation}

Next, let us consider $\int_{-\infty}^{0}e^{\gamma r}\alpha_{\tau, l}(r)dr$. The arguments we use are in part the same as those mentioned above. The function $r\mapsto \alpha_{\tau, l}(r)$ is continuous and uniformly  bounded with respect to both $\tau$ and $l$, thus the series $\sum_{n\geq 0}e^{-\gamma n}\alpha_{\tau,l}(r^*_n)$ is uniformly convergent, where $r^*_n=\displaystyle \argmax_{r \in [-n-1,-n]}\alpha_{\tau,l}(r)$. Therefore,
\begin{equation}\label{eq:estimate alpha_r}
\int_{-\infty}^{0}e^{\gamma r}\alpha_{\tau, l}(r)dr\leq \sum_{n\geq 0}e^{-\gamma n}\alpha_{\tau,l}(r^*_n). 
\end{equation}
Now, let $\epsilon>0$.  We can find an integer $N_{2}=N_{2}(\epsilon)>1$ such that
\begin{equation}\label{For: condmu}
\sum_{n> N_{2}}e^{-\gamma n}\alpha_{\tau,l}(r^*_n)<\frac{\epsilon^{p}}{2}.
\end{equation}
Furthermore, by putting $r^*=\displaystyle \argmax_{n \in \{0,\dots,N_2\}}\alpha_{\tau,l}(r^*_n)$, we have
\begin{eqnarray}
\sum_{n=0}^{ N_{2}}e^{-\gamma n}\alpha_{\tau,l}(r^*_n)&\leq &(1-e^{-\gamma})^{-1}\alpha_{\tau,l}(r^*).
\end{eqnarray}
Repeating the same argument as for $\alpha_{\tau,l}(0)$ enables us to write that for the positive number $l=l(\epsilon)$ provided by the condition $f\in \StAP[l]{p}_{K}(\R\times\espX, \espX)$ and for every $\tau$ belonging to the relatively dense set $\in T_{\St^{p}_{l}}(\epsilon,f,K)$ we have
\begin{equation}\label{eq:estimation_alpha_tau,l_general}
   \alpha_{\tau,l}(r^*) \leq (1-e^{-\gamma})\frac{\epsilon^{p}}{2}.
\end{equation}
Consequently, combining \eqref{eq:estimate alpha_r}, \eqref{For: condmu} and \eqref{eq:estimation_alpha_tau,l_general}, we obtain for the above $l$ and $\tau$ the desired estimation:
\begin{equation*}\label{For: equation alpha r}
\int_{-\infty}^{0}e^{\gamma r}\alpha_{\tau}(r)dr\leq \epsilon^{p}.
\end{equation*}
Which achieves the proof of Proposition~\ref{proposition: principale}.
}
\end{proof}
\medskip\medskip

{\em \textbf{Proof of Theorem \ref{theorem: principale}} We can break down the demonstration into two part. In the first part, we that, assuming
\eqref{For: L{Sp p>=2}}, equation \eqref{eq:EDO semi-lineaire} has a  unique  mild solution belonging to $\CB{(\R, \espX)}$, and given by \eqref{For: la solution de l'equation semi-linear}. In the second, using the technique by Kaminski et al. \cite{Kamenskii_Mikhail_and_Mellah_Omar_and_Raynaud_de_Fitte_15}, we show that solution \eqref{For: la solution de l'equation semi-linear} is Weyl almost periodic.
\subsection*{First part}
Clearly, following the same steps as in the proof of Theorem \ref{Theorem: linear}, $u$ is solution to \begin{equation*}
u(t) = \int_{-\infty}^{t}T(t-s)f(s,u(s))ds
\end{equation*}
if, and only if,
\begin{equation*}
u(t) = T(t-a)u(a) + \int_{a}^{t}T(t-s)f(s,u(s))ds,
\end{equation*}
for all $t\geq a$ for each $a\in\R$, which means that $u$ is a mild solution to \eqref{eq:EDO semi-lineaire}.

We introduce an operator $\Gamma$ by
\begin{equation*}
\Gamma u(t) = \int_{-\infty}^{t}T(t-s)f(s,u(s))ds.
\end{equation*}
Let us show that, assuming \eqref{For: L{Sp p>=2}}, the operator $\Gamma$ maps $\CB{(\R, \espX)}$ into $\CB{(\R, \espX)}$
and  has a unique  fixed point.

\textbf{First step.} Let us show that  $\Gamma u$ is bounded, for any $u\in \CB(\R, \espX)$.
From the hypotheses (H1), H\"older's inequality $(\frac{1}{p},
\frac{1}{q})$ and the triangular inequality, we have
\begin{align*}
\sup_{t\in \R}\norm{\Gamma u(t)} \leq& M\sup_{t\in
\R}\int_{-\infty}^{t}e^{-\delta(t-s)}\norm{f(s,u(s))}ds\\ \leq& M\sup_{t\in\R}
\sum_{k=1}^{\infty}\left(\int_{t-k}^{t-k+1}e^{-q\delta(t-s)}ds\right)^{\frac{1}{q}}\left(\int_{t-k}^{t-k+1}\norm{f(s,u(s))-
 f(s,0)}^{p}ds\right)^{\frac{1}{p}}\\& + M\sup_{t\in\R}
\sum_{k=1}^{\infty}\left(\int_{t-k}^{t-k+1}e^{-q\delta(t-s)}ds\right)^{\frac{1}{q}}
 \left(\int_{t-k}^{t-k+1}\norm{f(s,0)}^{p}ds\right)^{\frac{1}{p}}.
 \end{align*}
 Since the singleton $\{0\}$ is compact, $f(., 0)$ is bounded with respect to Stepanov's norm. By (H2) and H\"older's inequality, we have
 \begin{align*}
\sup_{t\in \R}\norm{\Gamma u(t)} \leq & M\sup_{t\in\R}\sum_{k=1}^{\infty}\left(\int_{t-k}^{t-k+1}e^{-q\delta(t-s)}ds\right)^{\frac{1}{q}}\left(\left(\int_{t-k}^{t-k+1}
\norm{L(s)}^{p}\norm{u(s)}^{p}ds\right)^{\frac{1}{p}} + \norm{f(., 0)}_{\St^{p}}\right)\\ \leq& M\frac{e^{\delta}}{e^{\delta}- 1} \sup_{t\in \R} \left(\left(\int_{t-k}^{t-k+1}\norm{L(s)}^{p}ds\right)^{\frac{1}{p}}\norm{u(.)}_{\infty} + \norm{f(., 0)}_{\St^{p}}\right)\\ \leq& M\frac{e^{\delta}}{e^{\delta}- 1}\Big(\norm{L(.)}_{\St^{p}}\norm{u(.)}_{\infty} + \norm{f(.,0)}_{\St^{p}}\Big) < \infty.
\end{align*}
Thus $\Gamma u$ is bounded.\medskip

\textbf{Second step} Let us show that the mapping $t \mapsto \Gamma u$ is continuous on $\R$, for every $u\in\CB(\R,\espX)$. For, let fix $u\in\CB(\R,\espX)$. We have
\begin{equation*}
\Gamma u(t) = \int_{-\infty}^{t}T(t-s)f(s,u(s))ds = \sum_{n=0}^{\infty}\int_{t-n-1}^{t-n}T(t-s)f(s,u(s))ds.
\end{equation*}
Let us set, for $n\in\N$,
\begin{equation*}
\Gamma_{n} u(t)=\int_{t-n-1}^{t-n}T(t-s)f(s,u(s))ds, \,t\in \R.
\end{equation*}
Let us check that $\Gamma_{n}(u)\in \CB(\R,\espX)$. To this end, let $t_{0}\in\R$.
Put $g(s)=: f(s, u(s))$. Clearly $g\in L_{\loc}^{p}(\R, \espX)$. Hence, in view of \cite[Teorem 2.4.2 p.70]{Alois}
\begin{equation}\label{eq:conti_g_mean_theo}
\lim_{h\rightarrow 0}\int_{t_{0}-n-1}^{t_{0}-n}\norm{g(s+h)- g(s)}^{p}ds = 0.
\end{equation}
Using H\"older's inequality, one has
\begin{align*}
\lim_{h\rightarrow 0}\|\Gamma_{n} u(t_{0}+h)&-\Gamma_{n} u(t_{0})\|\\ \leq&\lim_{h\rightarrow 0}M \int_{t_{0}-n-1}^{t_{0}-n}e^{-\delta (t_{0}-s)}\norm{g(s+h)- g(s)}ds\\ \leq& M\left(\int_{t_{0}-n-1}^{t_{0}-n}e^{-q\delta (t_{0}-s)}ds\right)^{\frac{1}{q}}\left(\lim_{h\rightarrow
0}\int_{t_{0}-n-1}^{t_{0}-n}\norm{g(s+h)- g(s)}^{p}ds\right)^{\frac{1}{p}}.
\end{align*}
The continuity of $\Gamma_{n} u$ is then a straightforward consequence of \eqref{eq:conti_g_mean_theo}.
The boundedness of $\Gamma_{n}u$ is a simple consequence of the following estimation
\begin{align*}
\norm{\Gamma_{n}u}_\infty \leq& \sup_{_{t\in\R}}\left(\int_{t-n-1}^{t-n}e^{-q\delta (t-s)}ds\right)^{\frac{1}{q}}\Big(\norm{L(.)}_{\St^{p}}\norm{u(.)}_{\infty} + \norm{f(.,0)}_{\St^{p}}\Big)\\ \leq & e^{-n\delta}\Big(\norm{L(.)}_{\St^{p}}\norm{u(.)}_{\infty} + \norm{f(.,0)}_{\St^{p}}\Big),
\end{align*}
from which we deduce in  addition that the series $\sum_{n=0}^{\infty}\norm{\Gamma_{n}u(t)}$ is uniformly convergent. The continuity of $\Gamma u$ is provided by the continuity of $\Gamma_{n}u$ for each $n$ and the uniform convergence of the series $\sum_{n\geq0}^{\infty}\norm{\Gamma_{n}u(t)}$.
\medskip
Thus $\Gamma$ maps $\CB(\R,\espX)$ into $\CB(\R,\espX)$.

\textbf{Third step} Let us show that, assuming \eqref{For: L{Sp p>=2}}, $\Gamma$ has a  fixed point in $\CB{(\R, \espX)}$. For every  $t\in \R$, we have from (H1)
\begin{align*}
\norm{\Gamma u(t) - \Gamma v(t)}&\leq
M\int_{-\infty}^{t}e^{-\delta(t-s)}\norm{f(s, u(s))- f(s,v(s))ds}\\
&\leq M\sum_{k =
1}^{+\infty}\int_{t-k}^{t-k+1}e^{-\delta(t-s)}\norm{f(s, u(s))-
f(s,v(s))}ds.
\end{align*}
Using H\"older's inequality and (H2), we get for every  $t\in \R$,
\begin{align*}
\norm{\Gamma u(t) - \Gamma v(t)} \leq&  M\sum_{k =
1}^{+\infty}\left(\int_{t-k}^{t-k+1}e^{-\delta q(t-s)}ds \right)^{\frac{1}{q}}\left(\int_{t-k}^{t-k+1}
\norm{L(s)}^{p}\norm{u(s)-v(s)}^{p}ds\right)^{\frac{1}{p}}\\
&\leq  M\sum_{k = 1}^{+\infty}\left(\int_{t-k}^{t-k+1}e^{-\delta
q(t-s)}ds\right)^{\frac{1}{q}}\left(\int_{t-k}^{t-k+1}\norm{L(s)}^{p}ds\right)^{\frac{1}{p}}\norm{u-v}_{\infty}\\
&\leq \left(\frac{M^{q}}{\delta
q}\right)^{\frac{1}{q}}\left(\frac{e^{\delta}}{e^{\delta}-1}\right)\norm{L}_{\St^{p}}\norm{u-v}_{\infty}.
\end{align*}
Hence, it follows that, for each $t\in\R$,
\begin{equation*}
\norm{\Gamma u - \Gamma v}_{\infty} \leq k\norm{u-v}_{\infty}.
\end{equation*}
Since  $k < 1$, we conclude that $\Gamma$ is a contraction operator.
We conclude that there exists a unique mild solution to \eqref{eq:EDO semi-lineaire} in $\CB{(\R, \espX)}$.
\subsection*{Second part}
In this part, we show that, under \eqref{For: L{Sp p>=2}} and \eqref{For: L{Sp p>2}}, the solution $u \in \CB{(\R, \espX)}$ to \eqref{eq:EDO semi-lineaire} is Weyl almost periodic.

By (H1), we obtain that, for any $\tau\in\R$ and any $t\in\R$,
\begin{align*}
\norm{u(t+\tau) - u(t)} &\leq \int_{-\infty}^{t}\norm{T(t-s)}\norm{f(s+ \tau,u(s+ \tau))-f(s,u(s))}ds\\ &\leq M\int_{-\infty}^{t}e^{-\delta(t-s)}\norm{f(s+ \tau,u(s+ \tau))-f(s,u(s))}ds.
\end{align*}
By  H\"older's inequality $(\frac{2}{p}, \frac{p-2}{p})$ and triangular inequality, we obtain
\begin{align*}
\norm{u(t+\tau) - u(t)} \leq& M\CCO{\frac{2p-4}{p\delta}}^{\frac{p-2}{p}} \CCO{\int_{-\infty}^{t}e^{-\frac{p\delta}{4}(t-s)}\norm{f(s+ \tau,u(s+ \tau))-f(s +\tau,u(s))}^{\frac{p}{2}}ds}^{\frac{2}{p}}\\&+ M\CCO{\frac{2p-4}{p\delta}}^{\frac{p-2}{p}} \CCO{\int_{-\infty}^{t}e^{-\frac{p\delta}{4}(t-s)}\norm{f(s+\tau, u(s))- f(s,u(t))}^{\frac{p}{2}}ds}^{\frac{2}{p}}\\
=: &M\CCO{\frac{2p-4}{p\delta}}^{\frac{p-2}{p}}\Big(I_{1}(t) + I_{2}(t)\Big).
\end{align*}
For $I_{1}(t)$, by (H2) and H\"older's inequality ($\frac{1}{2}, \frac{1}{2}$), we have
\begin{align*}
I_{1}(t) = &\CCO{\int_{-\infty}^{t}e^{-\frac{p\delta}{4}(t-s)}\norm{f(s+ \tau,u(s+ \tau))-f(s +\tau,u(s))}^{\frac{p}{2}}ds}^{\frac{2}{p}}\\
\leq& \CCO{\sum_{k = 0}^{+\infty}e^{-\frac{p\delta}{4}k}\int_{t-k-1}^{t-k} \norm{L(s + \tau)}^{p}ds}^{\frac{1}{p}}\CCO{\int_{-\infty}^{t}e^{-\frac{p\delta}{4}(t-s)}\norm{u(s + \tau)- u(s)}^{p}ds}^{\frac{1}{p}}\\ \leq&  \CCO{\frac{e^{\frac{p\delta}{4}}}{e^{\frac{p\delta}{4}} - 1}}^{\frac{1}{p}}\norm{L}_{\St^{p}}\CCO{\int_{-\infty}^{t}e^{-\frac{p\delta}{4}(t-s)}\norm{u(s +
\tau)- u(s)}^{p}ds}^{\frac{1}{p}}.
\end{align*}
Next, let us consider $I_{2}(t)$. By H\"older's inequality ($\frac{1}{2}, \frac{1}{2}$), we have
\begin{align*}
I_{2}(t) = &\CCO{\int_{-\infty}^{t}e^{-\frac{p\delta}{4}(t-s)}\norm{f(s+\tau, u(s))- f(s, u(s))}^{\frac{p}{2}}ds}^{\frac{2}{p}}\\ \leq&  \CCO{\int_{-\infty}^{t}e^{-\frac{p\delta}{4}(t-s)}ds}^{\frac{1}{p}} \CCO{\int_{-\infty}^{t}e^{-\frac{p\delta}{4}(t-s)}\norm{f(s+\tau, u(s))- f(s, u(s))}^{p}ds}^{\frac{1}{p}}\\  \leq& \CCO{\frac{4}{p\delta}}^{\frac{1}{p}} \CCO{\int_{-\infty}^{t}e^{-\frac{p\delta}{4}(t-s)}\norm{f(s+\tau, u(s))- f(s, u(s))}^{p}ds}^{\frac{1}{p}}.
\end{align*}
Thus, for all $\tau\in\R$, $\xi\in\R$ and $l>0$, we have
\begin{align*}
\frac{1}{l}\int_{\xi}^{\xi+l}\norm{u(t+\tau) -u(t)}^{p}dt \leq& M^{p}\CCO{\frac{2p-4}{p\delta}}^{p-2}\frac{1}{l}\int_{\xi}^{\xi+l}\norm{I_{1}(t) + I_{2}(t)}^{p}dt\\ \leq& M^{p}2^{p-1}\CCO{\frac{2p-4}{p\delta}}^{p-2}\CCO{\frac{1}{l}\int_{\xi}^{\xi+l} \norm{I_{1}(t)}^{p}dt+ \frac{1}{l}\int_{\xi}^{\xi+l} \norm{I_{2}(t)}^{p}dt}.
\end{align*}
Let us estimate the first term of the last inequality. We have, for every $\tau\in\R$, $\xi\in\R$ and $l>0$
\begin{align*}
\frac{1}{l}\int_{\xi}^{\xi+l} \norm{I_{1}(t)}^{p}dt \leq& \norm{L}_{\St^{p}}^{p} \CCO{\frac{e^{\frac{p\delta}{4}}}{e^{\frac{p\delta}{4}} - 1}}\frac{1}{l}\int_{\xi}^{\xi+l}\CCO{\CCO{\int_{-\infty}^{t}e^{-\frac{p\delta}{4}(t-s)}\norm{u(s + \tau)- u(s)}^{p}ds}^\frac{1}{p}}^{p}dt\\ \leq& \norm{L}_{\St^{p}}^{p}\CCO{\frac{e^{\frac{p\delta}{4}}}{e^{\frac{p\delta}{4}}-1}}\frac{1}{l}\int_{\xi}^{\xi+l}\int_{-\infty}^{t}e^{-\frac{p\delta}{4}(t-s)}
 \norm{u(s+\tau)- u(s)}^{p}ds\,dt.
\end{align*}
By a change of variables and Fubini's theorem, it follows that
\begin{equation}\label{For: equation I1}
\frac{1}{l}\int_{\xi}^{\xi+l} \norm{I_{1}(t)}^{p}dt \leq \norm{L}_{\St^{p}}^{p} \CCO{\frac{e^{\frac{p\delta}{4}}}{e^{\frac{p\delta}{4}} - 1}}\int_{-\infty}^{\xi}e^{-\frac{p\delta}{4}(\xi-s)}\CCO{\frac{1}{l}\int_{s}^{s+l}\norm{u(t + \tau)- u(t)}^{p}dt}ds.
\end{equation}
For the second term, using again a change of variables and Fubini's theorem, we have, for all $\tau\in\R$, $\xi\in\R$ and $l>0$,
\begin{align*}
\frac{1}{l}\int_{\xi}^{\xi+l} \norm{I_{2}(t)}^{p}dt \leq&  \CCO{\frac{4}{p\delta}} \frac{1}{l}\int_{\xi}^{\xi+l} \CCO{\CCO{\int_{-\infty}^{t}e^{-\frac{p\delta}{4}(t-s)}\norm{f(s+\tau, u(s))- f(s, u(s))}^{p}ds}^{\frac{1}{p}}}^{p}dt\\ \leq& \CCO{\frac{4}{p\delta}}\int_{-\infty}^{\xi}e^{-\frac{p\delta}{4}(\xi-s)} \frac{1}{l}\int_{s}^{s+l} \norm{f(t+\tau, u(t))- f(t, u(t))}^{p}dt\,ds\\=& \CCO{\frac{4}{p\delta}}\int_{-\infty}^{\xi}e^{-\frac{p\delta}{4}(\xi-s)}h_{\tau,l}(s)ds,
\end{align*}
where
\begin{equation*}
h_{\tau,l}(s):= \frac{1}{l}\int_{s}^{s+l} \norm{f(t+\tau, u(t))- f(t, u(t))}^{p}dt.
\end{equation*}
Summing up, for all $\xi, \tau\in\R$ and $l>0$, we obtain
\begin{equation*}
0\leq g_{\tau,l}(\xi):= \frac{1}{l}\int_{\xi}^{\xi+l}\norm{u(t+\tau) -u(t)}^{p}dt\leq C\alpha_{\tau,l}(\xi) + \beta \int_{-\infty}^{\xi}e^{-\delta_{1}(\xi-s)}g_{\tau,l}(s)ds,
\end{equation*}
with
\begin{equation*}
C:= \frac{M^{p}2^{p+1}}{p\delta}\CCO{\frac{2p-4}{p\delta}}^{p-2},\,\,\,\alpha_{\tau,l}(\xi):= \int_{-\infty}^{\xi}e^{-\frac{p\delta}{4}(\xi-s)} h_{\tau,l}(s)ds,\,\,\,\delta_{1}:= \frac{p\delta}{4},
\end{equation*}
and
\begin{equation*}
\beta:= M^{p}2^{p-1}\norm{L}_{\St^{p}}^{p} \CCO{\frac{2p-4}{p\delta}}^{p-2}\CCO{\frac{e^{\frac{p\delta}{4}}}{e^{\frac{p\delta}{4}} - 1}}.
\end{equation*}
It should be noted that the hypothesis \eqref{For: L{Sp p>2}} is equivalent $\delta_{1}>\beta$. We conclude by Lemma \ref{Lemma: Kame} that, for all $\gamma\in]0,\delta_{1}-\beta]$
such that $\int_{-\infty}^{0}e^{\gamma r}\alpha_{\tau,l}(r)dr$ converge, we have, for every $\xi\in\R$,
\begin{equation*}
0\leq g_{\tau,l}(\xi)\leq C\alpha_{\tau,l}(\xi) + C\beta \int_{-\infty}^{\xi}e^{-\gamma(\xi-r)}\alpha_{\tau,l}(r)dr.
\end{equation*}
We use Theorem~\ref{Theorem: Ursell} to get the Weyl almost periodicity of $u$. We take $\xi=0$, then for all $\tau\in\R$ and $l>0$, we have
\begin{equation*}
0\leq g_{\tau,l}(0)= \frac{1}{l}\int_{0}^{l}\norm{u(t+\tau) -u(t)}^{p}dt\leq C\alpha_{\tau,l}(0) + C\beta \int_{-\infty}^{0}e^{\gamma r}\alpha_{\tau,l}(r)dr.
\end{equation*}
 Therefore, in view of  Proposition \ref{proposition: principale}, for all $\epsilon>0$, there exist $l=l(\epsilon)>0$ and a relatively dense set $T_{\St^{p}_{l}}(\epsilon,f,K)$ such that
\begin{equation*}
\frac{1}{l}\int_{0}^{l}\norm{u(t+\tau) - u(t)}^{p}dt \leq C (1+ \beta)\epsilon^{p}, \quad \forall\tau\in T_{\St^{p}_{l}}(\epsilon,f,K).
\end{equation*}
Weyl almost periodicity of $u$ is then proved thanks to the inclusion $$T_{\St^{p}_{l}}\CCO{\epsilon^p C(1+ \beta),u}\supset T_{\St^{p}_{l}}\CCO{\epsilon^p,f,K}.$$

For $p=2$, the proof is similar to that for $p>2$. The details are omitted. \finpr
}
 {\em \remark {\em  The results obtained here can also be applied to the more general case $p \geq 1$ provided that $L (.)$ in (H1) is independent of $t$.
 }

\bibliographystyle{plain}{\em

}%

\end{document}